\documentclass[a4paper,12pt]{article}
\usepackage{times, url}
\usepackage{amsmath,amssymb,amsthm,amsfonts,bm,float,cite}

\newtheorem{theorem}{Theorem}[section]
\newtheorem{corollary}{Corollary}[section]
\newtheorem{definition}{Definition}[section]
\newtheorem{conjecture}{Conjecture}[section]
\newtheorem{remark}{Remark}[section]
\usepackage[none]{hyphenat}

\usepackage[center]{caption}

\begin{document}
\setcounter{page}{1}

\begin{center}
{\LARGE \bf  On some open problems concerning perfect powers\\[4mm] } 
\vspace{8mm}

{\Large \bf Marco Rip\`a}
\vspace{3mm}

sPIqr Society, World Intelligence Network\\ 
Rome, Italy\\
e-mail: \url{marcokrt1984@yahoo.it}
\vspace{2mm}

\end{center}
\vspace{10mm}

\noindent \sloppy {\bf Abstract:} The starting point of our paper is Kashihara's open problem number 30, concerning the sequence $A001292$ of the OEIS, asking how many terms are powers of integers. We confirm his last conjecture up to the $100128$-th term and provide a general theorem that rules out $4/9$ of the candidates. Moreover, we formulate a new, provocative, conjecture involving the OEIS sequence $A352991$ (which includes all the terms of $A001292$). Our risky conjecture states that all the perfect powers belonging to the sequence $A352991$ are perfect squares and they cannot be written as higher order perfect powers if the given term of $A352991$ is not equal to one. This challenging conjecture has been checked for any integer smaller than \linebreak
$10111121314151617181920212223456789$ and no counterexample has been found so far.
\vspace{5mm}

\section{Introduction} \label{sec:Intr}
In late 2010, the author of this paper found a recreative open problem by Kenichiro Kashihara (see \cite{kashihara:1}, open problem number 30, p. 25) concerning the sequence $A001292$ of the On-Line Encyclopedia of Integer Sequences (OEIS) \cite{A001292:2}. Kashihara's problem number 30 consists of two independent parts and the author solved the first one quite easily at the time (the complete solution can be found in \cite{ripa:6}, Section 3.3, pp. 12--15) since it asks to find the probability $0<p(c)<1$ that the trailing digit of the generic term of the sequence $A001292$ is $c\in \{0,1,2,\ldots ,9\}$ and the formula provided in \cite{ripa:6} shows that $p(c)=\frac{11-c}{55}$ for any $c\neq 0$, whereas $p(0)=0.0\overline{18}$ (e.g., if $c=7$, then $p(7)=\frac{4}{55} =0.0\overline{72}$).

In the present paper, we will focus on the second part of the above-mentioned Kashihara's problem number 30, asking how many elements of the sequence $A001292$ are perfect powers since Kashihara conjectured that there are none.

Now, bearing in mind that a perfect power of an integer $d>0$ is a natural number $k\geq{2}$ such that $a^k=d$, where also a is a positive integer, we could point out that $A001292(1)=1$ can be considered as a solution and argue how this disproves the conjecture, but (from here on) we will disregard this special case and assume that we are looking for a nontrivial counterexample to Kashihara's conjecture.

Lastly, Section 3 is devoted to introducing a (quite improbable) conjecture concerning perfect powers belonging to the OEIS sequence $A352991$ \cite{A352991:4, A353025:5}.

\section{The exclusion criterion}

To be clear on the invoked OEIS sequences, let us introduce a few useful definitions.

\begin{definition} \label{def1}
We define the $m$-th term of the OEIS sequence $A007908$ as \linebreak
$A007908(m) := 1\_2\_3\_\ldots\_(m-1)\_m$, where $m \in \mathbb{Z^+}$.
\end{definition}

\begin{definition} \label{def2}
We define the sequence $A001292$ of the OEIS as the concatenations (sorted in ascending order) of every cyclic permutation of the elements of the sequence $A007908$ \linebreak (e.g., given $m=3$, $A001292(A007908(3))=123,231,312$).
\end{definition}

\begin{definition} \label{def3}
We define the OEIS sequence $A352991$ as the concatenation of all the distinct permutations of the first strictly positive $m$ integers, sorted in ascending order (e.g., $12345671089$ is a term of the sequence, while $12345670189$ does not belong to $A352991$, even if all the digits of the string $1\_2\_3\_ \ldots \_9\_10$ appear once and only once, since ``$10$'' is missed).
\end{definition}

After having checked the first $100128$ terms of the sequence $A001292$ (see Appendix), exploring any exponent at or above two, we have not found any perfect power so that Kashihara's conjecture has been verified up to $10^{1235}$ (i.e., the $100129$-th term of $A001292$ is the smallest cyclic permutation of $A007908(448)$ and is greater than $10^{1235}$ by construction).
Moreover, we can prove the following Theorem \ref{theorem1}, concerning the sequence $A352991$ which includes every term of $A001292$.

\begin{theorem} \label{theorem1}
For any $m > 1$, $A352991(n)$ cannot be a perfect power of an integer if $A352991(n)$ is a permutation of $A007908(m)$ and $m : m \equiv \{2, 3, 5, 6\} \pmod 9$.
\end{theorem} 

\begin{proof} By definition, $A007908(m)$ \cite{A007908:3} cannot be a perfect power if $1\_2\_3\_\ldots \_(m-1)\_m$ is divisible by $3$ and it is not divisible by $3^2$. Thus, from the well-known divisibility by $3$ and $9$ criteria, $m : (3 \ | \sum_{j=1}^{m}  j) \wedge  (3^2 \nmid \sum_{j=1}^{m}  j) $ is a sufficient, but not necessary, condition for letting us disregard any permutation of $1\_2\_3\_\ldots \_(m-1)\_m$ (i.e., given $m$, if a generic permutation of $A007908(m)$ is divisible by $3$ and is not congruent to $0 \pmod 9$, then all the permutations of $A007908(m)$ are divisible by $3$ once and only once, since the commutativity property holds for addition).

It follows that, for any $n \in \mathbb{Z^+}$, $A352991(n)$ cannot be a perfect power if it is a permutation of the string $1\_2\_3\_\dots \_(m-1)\_m$, where m is such that $A134804(m)$ is divisible by $3$. Therefore, the residue modulo $9$ of every perfect power belonging to $A352991$ cannot be $2$ or $3$ or $5$ or $6$, and this concludes the proof of Theorem \ref{theorem1}. \end{proof}

\begin{corollary} \label{corollary1}
Kashihara's conjecture is true for the concatenation of any cyclic permutation of $A007908(m)$, where $m : m \equiv \{2, 3, 5, 6\}\pmod 9 \hspace{2mm} \vee \hspace{2mm} m < 448$.
\end{corollary}

\begin{proof} We observe that $A001292$ is a subsequence of $A352991$ \cite{A001292:2, A352991:4}. By invoking Theorem \ref{theorem1}, we can state that every perfect power candidate has to be the concatenation of a (cyclic) permutation of $A007908(m)$, where $m$ is such that $m : m \equiv \{2, 3, 5, 6\} \pmod 9$. On the other hand, all the remaining terms up to $99\_100\_101\_\dots \_445\_446\_447\_1\_2\_3\_\dots \_96\_97\_98$ have been directly checked (see Appendix for details) and no perfect power has been found.

Therefore, Corollary \ref{corollary1} confirms Kashihara's conjecture for any term of $A001292$ such that $m$ is congruent to $\{2, 3, 5, 6\} \pmod 9$ or $m \leq 447$.
\end{proof}

\begin{corollary} \label{corollary2}
 $\nexists n : A353025(n) \equiv \{2, 3, 4, 5, 6, 7, 8\} \pmod 9$, and any term of $A001292$ cannot be a perfect power if its digital root is not equal to $0$ or $1$.
\end{corollary}
 
\begin{proof}
Trivially, $10 \equiv 1 \pmod 9$ and also $(1+0) \equiv 1 \pmod 9$ so that any positive integer is congruent modulo $9$ to its digital root.

Now, we observe that every term of $A001292$ belongs to $A353025$.

Since from Theorem \ref{theorem1} it follows that every term of the sequence $A353025$ \cite{A353025:5} is a special permutation of $A007908(m)$ which is characterized by $m \equiv \{0, 1, 4, 7, 8\} \pmod 9$, in order to prove Corollary \ref{corollary2}, it is sufficient to note that
\begin{equation} \label{eq:1}
\sum_{j=1}^m j \equiv \begin{cases}
      0 \hspace{-2mm} \pmod 9  \quad   \text{ if } \quad m : m \equiv 0 \hspace{-2mm} \pmod  9\\
      1 \hspace{-2mm} \pmod  9  \quad  \text{ if }  \quad m : m \equiv 1 \hspace{-2mm} \pmod  9\\
      1 \hspace{-2mm} \pmod  9  \quad  \text{ if } \quad  m : m \equiv 4 \hspace{-2mm} \pmod  9\\
      1 \hspace{-2mm} \pmod  9  \quad  \text{ if } \quad  m : m \equiv 7  \hspace{-2mm} \pmod  9 \\
      0 \hspace{-2mm} \pmod  9  \quad  \text{ if } \quad  m : m \equiv 8 \hspace{-2mm} \pmod  9\\
    \end{cases}\,.
\end{equation}
\end{proof}

\begin{remark} \label{remark1} A well-known property of integers is that every perfect power that is congruent modulo $5$ to $0$ is also necessarily congruent to $\{0, 25, 75\}\pmod {100}$, while if a perfect power is congruent modulo $10$ to $6$, then its second last digit is odd.
\end{remark}

Thus, we are free to combine these additional constraints with Corollary \ref{corollary2} in order to reduce the number of perfect power candidates among the terms of $A352991$.

\section{Perfect cubes in A353025}
In the first half of April 2022, playing with Kashihara's conjecture, a more risky (very likely false but hard to disprove by brute force) conjecture arose, it is as follows.

\begin{conjecture} \label{conj1} Let $n \in \mathbb{N} - \text\{0, 1\}$ be given. We (provocatively) conjecture that if $n$ is such that $A352991(n)$ is a perfect power of an integer, then $\nexists k \in \mathbb{N} - \text\{0, 1, 2\} : A352991(n) = c^k, c \in \mathbb{N}$.
\end{conjecture}

On April 16, 2022, a direct search was performed by the author on the first $10^7$ terms of the sequence and no counterexample has been found ($42$ perfect squares only).
A few days later, Aldo Roberto Pessolano, performed a deeper search running the Mathematica codes published in the Appendix, without finding any counterexample and thus confirming Conjecture \ref{conj1} (at least) up to the smallest permutation of $A007908(22)$ (i.e., for any term of $A352991$ which is strictly greater than $1$ and smaller than $10111121314151617181920212223456789$).

\begin{remark} \label{remark2}
If confirmed, Conjecture \ref{conj1} would imply that all the perfect powers (strictly greater than $1$) in $A352991$ are perfect squares and nothing more (no cube, no square of square, and so forth). Nevertheless, under the (arbitrary, but perfectly reasonable) assumption of a standard probability distribution of the cubes in $A352991$ (i.e., we are assuming that $\frac{\mid n \in \mathbb{N} : A352991(n) \leq m  \wedge A352991(n)^{\frac{1}{3}} \mid}{\mid {n \in \mathbb{N} : A352991(n) \leq m}  \mid} \cong \frac{\mid n \in \mathbb{N} : A000578(n) \leq m \mid}{m}$ holds for any sufficiently large $m \in \mathbb{N}$), we would guess the existence of infinitely many counterexamples to Conjecture \ref{conj1}, even if the smallest one is expected to occur in the interval $[10^{58}, 10^{65}]$. On the other hand, the same argument would corroborate Kashihara's conjecture, since the number of perfect powers belonging to $A001292$ cannot probabilistically exceed
\begin{equation} \label{eq:3}
 2 \cdot \sum_{k=2}^{+\infty} \biggl(\sum_{j=309}^{+\infty}\frac{4 \cdot j  \cdot ( 10^{\frac{4 \cdot j }{k}} - 10^{\frac{4 \cdot j - 1}{k}})}{9 \cdot 10^{(4 \cdot j - 1)}}\biggr)
 \approx \frac{8}{9} \cdot \sum_{j=309}^{+\infty}\frac{ j  \cdot ( 10^{\frac{4 \cdot j}{2}} - 10^{\frac{4 \cdot j - 1}{2}})}{10^{(4 \cdot j - 1)}}
 \approx 0.
\end{equation}
\end{remark}

\noindent {\bf Additional open problems.} \textit{Does the sequence A353025 have infinitely many perfect squares, infinitely many perfect cubes, infinitely many perfect squares of squares, and so forth? Which is the smallest nontrivial perfect cube (if any) belonging to A353025 }(we point out that all the terms greater than one and below $(1.01\cdot10)^{40}$ have been checked without finding any cube)\textit{?}

\section{Conclusion}
Kashihara's open problem number 30 has not been completely solved yet. Even if the first part, concerning the probability that the trailing digit of $A001292(n)$ is $c=1,2,\dots ,9$, was solved by the author a dozen years ago \cite{ripa:6}, the conjecture in the second part still need a proof or a nontrivial counterexample (the smallest candidate has $1236$ digits).
Moreover, in the present paper, we have introduced a wider speculation that allows us to ask ourselves how to find a term of the OEIS sequence $A353025$ (disregarding $A353025(1))$ which is not a perfect square; a challenging open problem, considering that there is not any perfect cube among the terms of $A352991$ in the interval $(1,10^{40}]$.

\section{Appendix}

Aldo Roberto Pessolano helped the author of the present paper by verifying Kashihara's conjecture and Conjecture \ref{conj1} for a very large number of terms. All the provided Mathematica codes run on the M1 processor of his Apple MacBook Air (2020).
Kashihara's conjecture has been currently tested up to the $100128$-th term of $A001292$ and we confirm that it holds for every element of the set ${A001292(2), A001292(3),\ldots, A001292(100128)}$.

The search reached the term $99\_100\_\dots \_446\_447\_1\_2\_\ldots \_97\_98 \approx 9.91 \cdot 10^{1232}$ in $28823$ seconds (about $8$ hours of calculations) and the code is as follows:
\begin{verbatim}
c = True;
p = Table[Prime[q], {q, 1, 565}];
Do[rn = Range[k];
n = ToExpression[StringJoin[ToString[#]&/@rn]];
If[And[Mod[n, 9] != 3, Mod[n, 9] != 6],
Do[r = RotateLeft[rn, i - 1];
nk = ToExpression[StringJoin[ToString[#]&/@r]];
If[IntegerQ[nk^(1/#)],
Print[nk, " = ", nk^(1/#), "^", #]; c = False; Break[]
]&/@p,
{i, 1, k}]
];
If[c, Print["1..", k, " checked."], Break[]],
{k, 2, 447}]
\end{verbatim}
\vspace*{2mm}

About our investigation on the perfect powers in $A352991$, Pessolano has recently completed the direct check of every term of $A352991$ which falls in the interval
$(1, 987654322120191817161514131211110]$ (see the code below).

As expected, the test has not returned any perfect power above two.
\begin{verbatim}
z = False;
h = 3;
p = Table[Prime[q], {q, 2, 10}];
q[x_, k_, d_, m_] :=(y = x^k;
	If[DigitCount[y] == d,c = True;
		Do[If[Not[StringContainsQ[ToString[x], ToString[i]]],c = False;
			Break[],c = True],{i, 10, m}],c = False];
		Return[c])
	Do[r = Range[k];
		n = ToExpression[StringJoin[ToString[#]&/@r]];
		If[And[Mod[n, 9] != 3, Mod[n, 9] != 6],d = DigitCount[n];
			(s = IntegerPart[(10^(IntegerLength[n] - 1))^(1/#)];
				f = IntegerPart[(10^(IntegerLength[n]))^(1/#)];
				Do[If[q[x, #, d, k], Print[x, "^", #, " = ", y];
					z = True; Break[]],{x, s, f}])&/@p;g = 2^h;
	While[g < n,If[q[#, h, d, k], 
		Print[x, "^", h, " = ", y]; 
		z = True; 
		Break[]]&/@{2, 3, 5, 6, 7};
	h++;g = 2^h]];
If[z, Break[], 
	Print["1..", k, " checked."]],
{k, 2, 21}]

\end{verbatim}

On the other hand, the following code returns the complete list of the smallest $42$ perfect squares belonging to $A352991$.
\begin{verbatim}
z = 1;
Do[r = Range[k];
n = ToExpression[StringJoin[ToString[#]&/@r]];
If[And[Mod[n, 9] != 3, Mod[n, 9] != 6],
d = DigitCount[n];
s = IntegerPart[Sqrt[10^(IntegerLength[n] - 1)]];
f = IntegerPart[Sqrt[10^(IntegerLength[n])]];
Do[y = x^2;
If[DigitCount[y] == d,
c = True;
Do[
If[Not[StringContainsQ[ToString[y], ToString[i]]],
c = False
],
{i, 10, k}];
If[c, Print[z, " ", y]; z++]
],
{x, s, f}]
],
{k, 2, 10}]

\end{verbatim}

These $42$ perfect squares correspond to all the perfect powers in
$(1, 10^{16}]$ belonging to $A352991$, while the next perfect square is $10135681742311129$ (we observe that $100676123^2$ is a permutation of $1\_2\_3\_\ldots \_16$, as suggested by the statement of Theorem \ref{theorem1}). \vspace{2mm}

\noindent1 	\textbf{    }	$13527684$\\
2	\textbf{    }	$34857216$\\
3 	\textbf{    }	$65318724$\\
4 	\textbf{    }	$73256481$\\
5 	\textbf{    }	$81432576$\\
6 	\textbf{    }	$139854276$\\
7 	\textbf{    }	$152843769$\\
8 	\textbf{    }	$157326849$\\
9 	\textbf{    }	$215384976$\\
10 	\textbf{    }	$245893761$\\
11 	\textbf{    }	$254817369$\\
12 	\textbf{    }	$326597184$\\
13 	\textbf{    }	$361874529$\\
14 	\textbf{    }	$375468129$\\
15 	\textbf{    }	$382945761$\\
16 	\textbf{    }	$385297641$\\
17 	\textbf{    }	$412739856$\\
18	\textbf{    }	$523814769$\\
19 	\textbf{    }	$529874361$\\
20 	\textbf{    }	$537219684$\\
21 	\textbf{    }	$549386721$\\
22 	\textbf{    }	$587432169$\\
23 	\textbf{    }	$589324176$\\
24 	\textbf{    }	$597362481$\\
25 	\textbf{    }	$615387249$\\
26 	\textbf{    }	$627953481$\\
27 	\textbf{    }	$653927184$\\
28 	\textbf{    }	$672935481$\\
29 	\textbf{    }	$697435281$\\
30 	\textbf{    }	$735982641$\\
32 	\textbf{    }	$743816529$\\
33 	\textbf{    }	$842973156$\\
34 	\textbf{    }	$847159236$\\
35 	\textbf{    }	$923187456$\\
36 	\textbf{    }	$14102987536$\\
37 	\textbf{    }	$24891057361$\\
38 	\textbf{    }	$27911048356$\\
39 	\textbf{    }	$28710591364$\\
40 	\textbf{    }	$57926381041$\\
41 	\textbf{    }	$59710832164$\\
42 	\textbf{    }	$75910168324$
\newline

In the end, our tests have finally confirmed that all the perfect powers that are smaller than $10^{34}$ and that belong to the OEIS sequence $A352991$ are perfect squares (only). At present, Conjecture \ref{conj1} has been tested for every integer smaller than $10111121314151617181920212223456789$, and no counterexample has been found yet.

\section*{Acknowledgments} 

We sincerely thank Aldo Roberto Pessolano for having helped us very much with the search for perfect powers belonging to $A001292$ and $A352991$, letting us confirm Kashihara's conjecture up to $10^{1235}$ and Conjecture \ref{conj1} up to $(1.01\cdot10)^{34}$.
 
\makeatletter
\renewcommand{\@biblabel}[1]{[#1]\hfill}
\makeatother

\bibliographystyle{plain}
\bibliography{On_some_open_problems_concerning_perfect_powers}

\end{document}